\documentclass[12pt]{article}

\usepackage{amsfonts}
\usepackage{amsmath}
\usepackage{tikz}
\usepackage{graphicx}

\makeatletter
\newcommand*{\encircled}[1]{\relax\ifmmode\mathpalette\@encircled@math{#1}\else\@encircled{#1}\fi}
\newcommand*{\@encircled@math}[2]{\@encircled{$\m@th#1#2$}}
\newcommand*{\@encircled}[1]{%
  \tikz[baseline,anchor=base]{\node[draw,circle,outer sep=0pt,inner sep=.2ex] {#1};}}
\makeatother

\begin{document}
\title{Chebyshev polynomials in the 16th century}
\author{Walter Van Assche\footnote{Supported by Research Foundation Flanders (FWO), Grant G.0C9819N.} \\ KU Leuven, Belgium}
\date{March 18, 2022}
\maketitle

\begin{abstract}
We give a few examples of Chebyshev polynomials that appeared in mathematical problems from the 16th and 17th century.
The main example is the famous equation of Adrianus Romanus (Adriaan van Roomen) containing a polynomial of degree $45$.
\end{abstract}

\section{Introduction}
\begin{figure}[h!]
\includegraphics[width=6in]{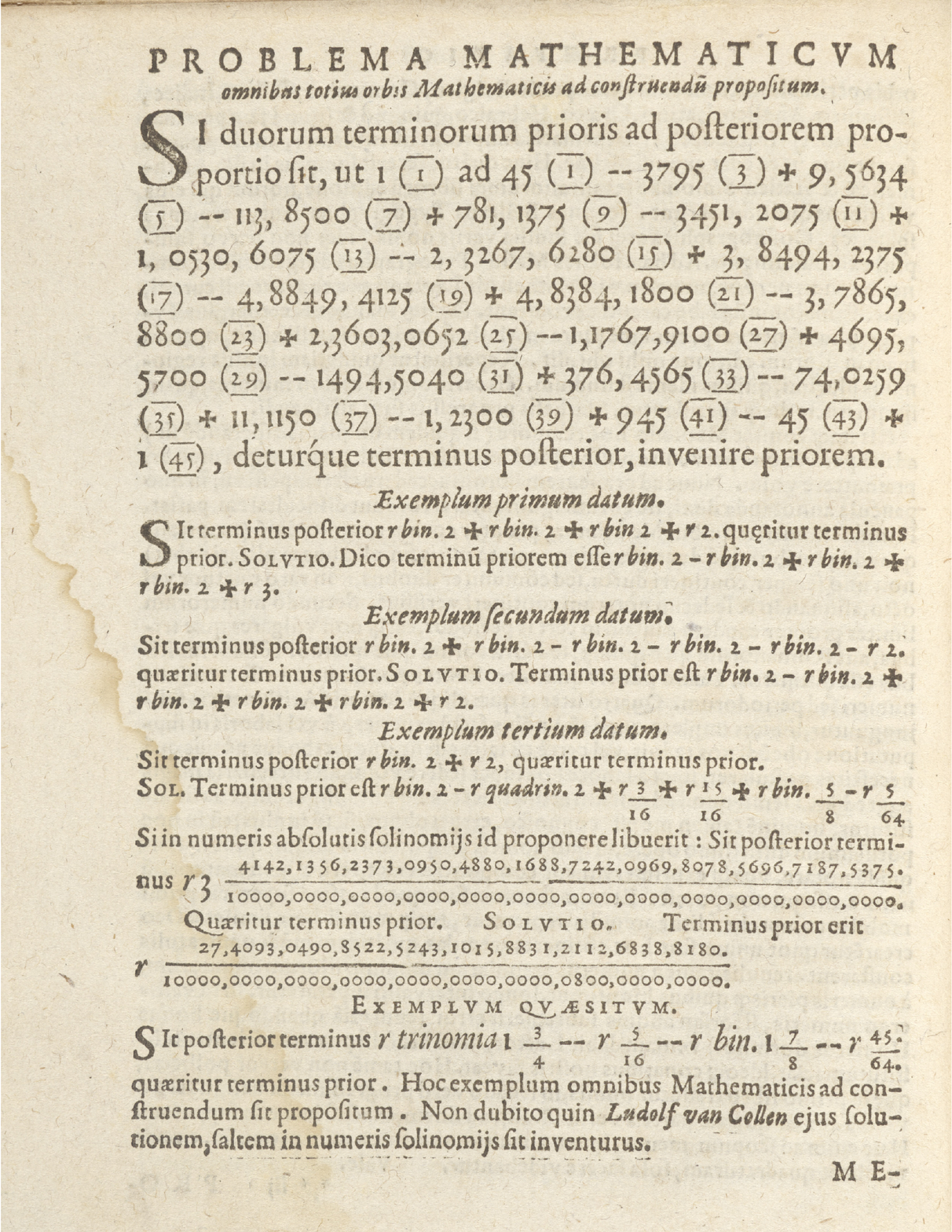}
\caption{Page from Adrianus Romanus \cite{Romanus}, KU Leuven BRES: Tabularium, CaaA713}
\label{Romanus-fig1}
\end{figure}

\begin{figure}[h!]
\centering
\includegraphics[width=4in]{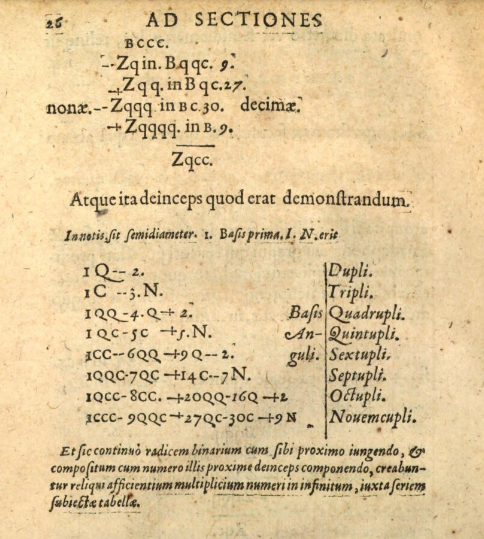}
\caption{Page 26 from Vi\`ete's \textit{Ad angulares sectiones} \cite{Viete1615}.}  Bayerische Staatsbibliothek M\"unchen, 4 A.gr.b 128\#Beibd.2 \\
(Digital Object ID bsb10215489).
\label{figV1}
\end{figure}

\begin{figure}[ht]
\centering
\includegraphics[width=4in]{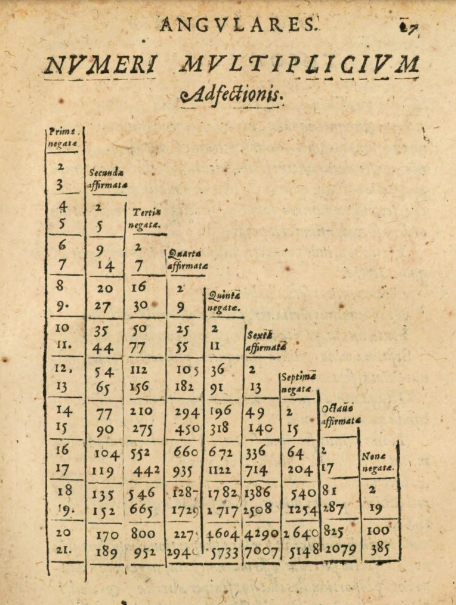}
\caption{Page 27 from Vi\`ete's \textit{Ad angulares sectiones} \cite{Viete1615}.}  Bayerische Staatsbibliothek M\"unchen, 4 A.gr.b 128\#Beibd.2 \\
(Digital Object ID bsb10215489).
\label{figV2}
\end{figure}

Recently I came accross a paper \cite{Bockstaele} (in Dutch) by Paul Bockstaele (1920--2009), who was my professor for infinitesimal analysis (calculus) and differential
equations when I was a student in the late 1970's. Bockstaele was an expert in the history of mathematics and a specialist in the work of Adrianus
Romanus (Adriaan van Roomen, Leuven 1561 -- Mainz 1615), a mathematician and professor of medicin at the University of Louvain (1586--1592) and the University of W\"urzburg (1593--1603). 
He was famous for computing the value of $\pi$ to sixteen decimals in 1593, which was better than the ten digits that Fran\c{c}ois Vi\`ete had obtained earlier. 
Romanus' record was later beaten by his friend Ludolf van Ceulen (1540--1610), who calculated 20 decimals in 1596 and later reached 35 decimals, 
which were engraved on his tombstone in the Pieterskerk in Leiden\footnote{A replica containing this approximation of $\pi$ was placed in the Pieterskerk in 2000, see e.g. \cite{Pieterskerk}.}.
Bockstaele's paper drew my attention because it contained a polynomial of degree 45 from a problem posed by Adrianus Romanus in 1593, see
Fig. \ref{Romanus-fig1}. In order to recognize a polynomial, one needs to know that Romanus was using notation introduced earlier by Simon Stevin,
where the unknown $x$ is denoted by a circle $\encircled{1}$ with the number 1 inside, and powers of $x$ were denoted by a circle containing the power, e.g. $\encircled{7}$ stands for
$x^7$ and $\encircled{45}$ stands for $x^{45}$ \cite[p. 157]{Cajori}. Stevin's notation is difficult for printing and hence Romanus uses instead two rounded
parentheses with horizontal lines above and below \cite[p.~344]{Cajori}.
The polynomial in Fig.  \ref{Romanus-fig1} is therefore
\begin{align*}
 &45 x - 3795 x^3 + 95634 x^5 - 1138500 x^7 + 7811375 x^9 - 34512075 x^{11} + 105306075 x^{13} \\
 &- 232676280 x^{15} + 384942375 x^{17}    - 488494125 x^{19} + 483841800 x^{21} - 378658800 x^{23} \\
 &+ 236030652 x^{25} - 117679100 x^{27} + 46955700 x^{29} - 14945040 x^{31}     + 3764565 x^{33} \\
 &- 740259 x^{35} + 111150 x^{37} - 12300 x^{39} + 945 x^{41} - 45 x^{43} + x^{45} .  
\end{align*}    
\newpage
\newpage
This polynomial has integer coefficients and contains only odd powers, so that it corresponds to an odd function. 
I recognized this as a Chebyshev polynomial of the first kind with a scaling of the variable, in particular this is $2 T_{45}(x/2)$.
Of course I am not the first to make this observation: nowadays we know that the Chebyshev polynomials are related to the trigonometric functions by
\begin{equation}  \label{Tn-cos}
     T_n(x) = \cos n\theta, \qquad x = \cos \theta, 
\end{equation}     
so the Chebyshev polynomial of the first kind expresses how to obtain $\cos n\theta$ from $\cos \theta$. 
For odd values of the degree, say $n=2m+1$, one finds
\begin{equation}  \label{Tn-sin}
   T_{2m+1}(\sin \theta) = (-1)^m \sin (2m+1)\theta  
\end{equation}   
and this formula makes the connection between Chebyshev polynomials and regular polygons, because the side of a regular $n$-gon is $2\sin (\pi/n)$.
In particular, if one puts $\theta=\pi/(2m+1)$, then \eqref{Tn-sin} shows that $2\sin(\pi/(2m+1))$ is a root of the equation $2T_{2m+1}(x/2)=0$ and it is the
smallest positive root. Such expressions were already found by
Fran\c{c}ois Vi\`ete (1540--1603) but his work on calculating $\cos n\theta$ and $\sin n\theta$ in terms of $\cos \theta$ and $\sin \theta$ was published
only after his death in 1615 \cite{Viete1615}, see  Figs. \ref{figV1} and \ref{figV2}, where another notation is used for powers of the variable $x$, 
with  $ N = x$, $Q = x^2$, $C = x^3$, $QQ = x^4$, $QC=x^5$, $CC = x^6$, etc.
In Fig. \ref{figV1} one recognizes $2T_n(x/2)$ for $n=2,3,\ldots,9$, and in Fig. \ref{figV2} the coefficients are given up to $n=21$ (with a few misprints).

\noindent  Naturally he did not use the terminology \textit{Chebyshev polynomials}, because Pafnuty L. Chebyshev (1821--1894) only appears two centuries later. 
 Vi\`ete recognized the polynomial of degree $45$ to be the one expressing $2 \sin (45\theta)$ in terms of $2\sin \theta$, as in \eqref{Tn-sin}.

\section{The problem of Adrianus Romanus}  \label{sec2}
Adrianus Romanus spent quite a large part of his life to work out a method to compute the sides, circumference and area of a regular polygon
inscribed in a circle of radius one. The purpose was to approximate the circumference of the circle (which is $2\pi$) by the circumference of a regular polygon with $n$ sides,
by taking $n$ sufficiently large. This idea started with Archimedes and was also the method used by Vi\`ete and van Ceulen for their approximation of $\pi$.
The problem mainly was to compute the length of one side of a regular $n$-gon for large $n$. Such a side has length $2 \sin (\pi/n)$, and the circumference of the $n$-gon is $2n \sin (\pi/n)$. The limit
\[   \lim_{n \to \infty} 2n \sin (\pi/n) = 2\pi  \]
indeed shows that the circumference of the regular $n$-gon tends to the circumference of the circle.
Adrianus Romanus intended to publish a series of books devoted to this problem, but only four were published. At the end of the preface (pr\ae fatio) of his first book \cite{Romanus}, he
proposed a \textit{Mathematical problem for all the mathematicians in the world}, see Fig. \ref{Romanus-fig1}. In contemporary notation, the problem is to solve the equation
\begin{equation}   \label{T45=b}
      2T_{45}(x/2) = b  
\end{equation}
when $b$ is a given expression. He gave three examples with a solution. Again one has to get used to the notation used in the original document: a square root was denoted
by $\surd$ (which looks a bit like an $r$), and $\surd\, bin[omia]$ indicates a square root with two terms, $\surd\,trinomia$ indicates a square root with three terms, $\surd\,quadrin[omia]$ indicates a square root
with four terms. The three examples are 
\begin{description}
\item[First example:] when $b=\sqrt{2+\sqrt{2+\sqrt{2 + \sqrt{2}}}}$, then the solution is 
\[    x = \sqrt{2 - \sqrt{2 + \sqrt{2 + \sqrt{2 + \sqrt{3}}}}}.  \]
\item[Second example:] when $b=\sqrt{2+\sqrt{2-\sqrt{2-\sqrt{2-\sqrt{2-\sqrt{2}}}}}}$, then the solution is
\[      x = \sqrt{2-\sqrt{2+\sqrt{2+\sqrt{2+\sqrt{2+\sqrt{3}}}}}}.  \]
\item[Third example:] when $b = \sqrt{2+\sqrt{2}}$, then the solution is
\[   x = \sqrt{2- \sqrt{2 + \sqrt{\frac{3}{16}} + \sqrt{\frac{15}{16}} + \sqrt{\frac58 - \sqrt{\frac{5}{64}}}}}. \]
\end{description}
The actual problem was to solve the equation  \eqref{T45=b} when $b= \sqrt{1\frac34 - \sqrt{\frac{5}{16}} - \sqrt{1\frac78 - \sqrt{\frac{45}{64}}}}$.
At first sight this looks like an impossible problem since it requires to solve a polynomial equation with a polynomial of degree 45, and we only know how to solve
polynomial equations of degree two, three or four. However, for some polynomial equations of degree $\geq 5$ one can still find solutions and in this particular case
the polynomial is related to trigonometric functions and one may use a geometric reasoning to find a solution. For instance, one can solve the equation
$T_n(x)=0$ for a Chebyshev polynomial of the first kind for any $n$, and the solution is
\[     x = \cos \left( \frac{(2k-1)\pi}{2n}\right), \qquad    k=1,2,\ldots,n,  \]
which for $n=2^j$ can be expressed in radicals using $j$ nested square roots (see Section \ref{sec5}).
So one really needs to take into account that \eqref{T45=b} contains a Chebyshev polynomial and a special choice of $b$.
Another peculiar thing is that Adrianus Romanus seems to be looking for one solution, but the equation \eqref{T45=b} has $45$ solutions due to the fundamental theorem of algebra.
The latter theorem is a famous result of Gauss, but he appeared only 200 years after Romanus. If $0 < b < 2$ then there are indeed 45 real solutions to the equation \eqref{T45=b},
and Romanus was only interested in the smallest positive solution, which corresponds to the length of one side of some regular polygon, depending
on the choice of $b$.

Fran\c{c}ois Vi\`ete found the solution shortly after its announcement and published it in 1595 \cite{Viete1595}. 
In fact, Vi\`ete gave all the positive solutions of the equation (there are 23 positive
and 22 negative solutions). Vi\`ete also noted that the second example is wrong and suggested a correction: he proposed to take
$b=\sqrt{2-\sqrt{2-\sqrt{2+\sqrt{2+\sqrt{2+\sqrt{2}}}}}}$, but this also does not give the solution given by Romanus.

\section{The solution of the problem}  \label{sec3}
Before we work out the solution of the actual problem, we will have a look at the three examples. Recall that for some special values of $\theta$ one can find
$\cos \theta$ and $\sin \theta$ using only square roots. 

\begin{table}[ht]
\centering
\begin{tabular}{|c|c|c|}

\hline
$n$ & $\cos \frac{2\pi}{n}$ & $\sin \frac{2\pi}{n}$  \rule{0pt}{15pt} \\[12pt]
\hline
$1$ &  $1$         & $0$   \rule{0pt}{15pt}\\
$2$ & $-1$         & $0$   \rule{0pt}{15pt} \\
$3$ & $-\frac12$ & $\frac{\sqrt{3}}{2}$  \rule{0pt}{15pt} \\
$4$ & $0$          & $1$ \rule{0pt}{15pt} \\
$5$ & $\frac{\sqrt{5}-1}{4}$ & $\frac{\sqrt{10+2\sqrt{5}}}{4}$ \rule{0pt}{15pt}\\
$6$ & $\frac12$  & $\frac{\sqrt{3}}{2}$  \rule{0pt}{15pt} \\
$7$  & -- & --  \rule{0pt}{15pt} \\
$8$ & $\frac{\sqrt{2}}{2}$ & $\frac{\sqrt{2}}{2}$  \rule{0pt}{15pt} \\
$9$ & -- & --  \rule{0pt}{15pt}  \\
$10$ & $\frac{\sqrt{5}+1}{4}$ & $\frac{\sqrt{10-2\sqrt{5}}}{4}$  \rule{0pt}{15pt} \\
$11$ & -- & -- \rule{0pt}{15pt} \\
$12$ & $\frac{\sqrt{3}}{2}$ & $\frac12$ \rule{0pt}{15pt} \\
\hline
\end{tabular}
\caption{Values for special angles}
\label{special-angles} 
\end{table}

For instance, for $\sin (2\pi/5)$ one can use \eqref{Tn-sin} to find
\[    T_5(\sin \theta) = \sin(5 \theta), \]
and using $T_5(x) =16x^5 - 20x^3 + 5x$ and $\sin 2\pi = 0$, one then finds
\[        16 \sin^5 (2\pi/5) - 20 \sin^3 (2\pi/5)  +5 \sin(2\pi/5) = 0. \] 
This quintic equation has 5 solutions: $0$, two negative and two positive roots and $\sin (2\pi/5)$ corresponds with the largest positive solution
\[   \sin \frac{2\pi}{5} = \frac{\sqrt{10+2\sqrt{5}}}{4} .  \]
The value for $\cos (2\pi/5)$ can then be found
from $\cos^2\theta + \sin^2\theta = 1$. For $n=7$ and $n=9$  the result will involve
third roots.
One can then use these special angles and simple trigonometric formulas to find the value of $\cos (2\pi/n)$ and $\sin (2\pi/n)$ for other
values of $n$. Recall
\begin{equation} \label{cos2}
     2 \cos (\theta/2) = \sqrt{2 + 2 \cos \theta},  
\end{equation}
and
\begin{equation}  \label{sin2}
    2  \sin (\theta/2) = \sqrt{2 - 2 \cos \theta},     
\end{equation}
then starting with $2\cos (\pi/4) = \sqrt{2}$ one finds from \eqref{cos2} that
\[    2 \cos (\pi/8) = \sqrt{2+\sqrt{2}}, \]
which is the value of $b$ in Romanus' third example. Using \eqref{cos2} two more times gives
\[    2 \cos (\pi/32) = \sqrt{2 + \sqrt{2 + \sqrt{2 + \sqrt{2}}}} , \]
which is the value of $b$ in Romanus's first example.

If we use \eqref{Tn-sin}, then the equation \eqref{T45=b} for the first example becomes
\[     2 \sin (45 \theta) = 2 \cos \frac{\pi}{32} = 2 \sin \frac{15\pi}{32}, \qquad   x = 2 \sin \theta,  \]
and the solution is readily seen to be
\[         \theta =  \frac{15\pi}{45 \times 32} + \frac{2\pi k}{45}, \qquad   k \in \mathbb{Z} .  \]
The solutions for \eqref{T45=b} are therefore given by
\[     x = 2 \sin \left( \frac{\pi}{96} + \frac{2\pi k}{45} \right), \qquad k \in \mathbb{Z}.  \]
This gives 22 negative solutions and 23 positive solutions. The smallest positive solution is
\[   x = 2 \sin \frac{\pi}{96}, \]
which corresponds to the length of one side of a regular polygon with $96$ sides.
Now, if we start from $2\cos( \pi/6) = \sqrt{3}$ and apply \eqref{cos2} three times, then
\[   2 \cos \frac{\pi}{48} = \sqrt{2+\sqrt{2+\sqrt{2+\sqrt{3}}}}, \]
and then applying \eqref{sin2} gives
\[    2 \sin \frac{\pi}{96} = \sqrt{2-\sqrt{2+\sqrt{2+\sqrt{2+\sqrt{3}}}}}, \]
which is indeed the expression that Adrianus Romanus gave as the solution.

For the third example the equation \eqref{T45=b} becomes
\[   2 \sin (45\theta) = 2 \cos (\pi/8) = 2 \sin \frac{3\pi}{8}, \qquad   x = 2 \sin \theta, \]
and hence 
\[   \theta = \frac{3\pi}{45 \times 8} + \frac{2\pi k}{45}, \qquad k \in \mathbb{Z}. \]
The solutions are then given by
\[   x = 2 \sin \left(  \frac{\pi}{120} + \frac{2\pi k}{45} \right), \qquad k \in \mathbb{Z}, \]
which gives 45 solutions, 22 negative and 23 positive solutions, and the smallest positive solution is
\[   x =  2 \sin \frac{\pi}{120}, \]
which corresponds to the length of one side of a regular polygon with 120 sides.
A simple exercise gives
\[   \cos \frac{\pi}{30} = \cos \left( \frac{\pi}{5}-\frac{\pi}{6} \right) = \cos \frac{\pi}{5} \cos \frac{\pi}{6} + \sin \frac{\pi}{5} \sin \frac{\pi}{6}, \]
so that (see Table \ref{special-angles})
\[    2 \cos \frac{\pi}{30} = \frac{\sqrt{3}(1+\sqrt{5})}{4} + \frac{\sqrt{10-2\sqrt{5}}}{4} 
            = \sqrt{\frac{3}{16}} + \sqrt{\frac{15}{16}} + \sqrt{\frac58 - \sqrt{\frac{5}{64}}} .  \] 
If we then use \eqref{cos2} followed by \eqref{sin2}, we find that
\[   2 \sin \frac{\pi}{120} = \sqrt{2 - \sqrt{2+   \sqrt{\frac{3}{16}} + \sqrt{\frac{15}{16}} + \sqrt{\frac58 - \sqrt{\frac{5}{64}}}}}, \]
which again corresponds to the expression for the solution of example three given by Adrianus Romanus.   

So we now see the pattern and we can go for the actual problem for which
\[  b = \sqrt{\frac{7}{4} - \sqrt{\frac{5}{16}} - \sqrt{\frac{15}{8}-\sqrt{\frac{45}{64}}}} = \sqrt{2 - \frac{1+\sqrt{5}}{4} - \sqrt{3}\frac{\sqrt{10-2\sqrt{5}}}{4}} . \]                 
First we observe that
\[    \cos \frac{2\pi}{15} = \cos \left( \frac{\pi}{3}-\frac{\pi}{5} \right) = \cos \frac{\pi}{3} \cos \frac{\pi}{5} + \sin \frac{\pi}{3} \sin \frac{\pi}{5}, \]
so that Table \ref{special-angles} gives
\[    2 \cos \frac{2\pi}{15} = \frac{1+\sqrt{5}}{4} + \sqrt{3}\frac{\sqrt{10-2\sqrt{5}}}{4}. \]
Then use \eqref{sin2} to find that
\[   b = 2 \sin \frac{\pi}{15} .  \]
Equation \eqref{T45=b} then becomes
\[    2 \sin (45\theta) = 2 \sin \frac{\pi}{15}, \qquad   x = 2\sin \theta, \]
and therefore
\[    \theta = \frac{\pi }{15\times 45} + \frac{2\pi k}{45}, \qquad k \in \mathbb{Z}, \]
which gives the solutions
\[    x = 2 \sin\left( \frac{\pi }{675} + \frac{2\pi k}{45} \right), \qquad k \in \mathbb{Z}.  \]
The smallest positive solution is
\[    x =  2 \sin \frac{\pi}{675},  \]
which corresponds to the length of one side of a regular polygon with $675$ sides.
This number cannot be written in terms of square roots only because $675 = 3^3 5^2$ and hence one needs solving cubic and quintic equations (see Section \ref{sec5}).

The second example contains an error, as was observed by Vi\`ete. The solution given by Adrianus Romanus is
\[  x =  2 \sin \frac{\pi}{192}, \]
which corresponds to the length of one side of a regular polygon with 192 sides. But if we compute $2T_{45}(x/2)$ for this value of $x$, then we find
\[   2T_{45}\left( \sin \frac{\pi}{192} \right) = 1.3431179096940368013\ldots \]
and the value for $b$ given by Romanus is
\[   b = \sqrt{2+\sqrt{2-\sqrt{2-\sqrt{2-\sqrt{2-\sqrt{2}}}}}} = 1.7401739822174228373\ldots \]
The value for $b$ proposed by Vi\`ete is also not correct:
\[   b= \sqrt{2-\sqrt{2-\sqrt{2+\sqrt{2+\sqrt{2+\sqrt{2}}}}}} = 1.3790810894741338492\ldots  \]
A correct value for $b$ is
\[    b = \sqrt{2-\sqrt{2-\sqrt{2+\sqrt{2+\sqrt{2}}}}} ,  \]
which corresponds to $b=2\sin (15\pi/64)$.
\newpage
\begin{figure}[h!]
\includegraphics[width=6in]{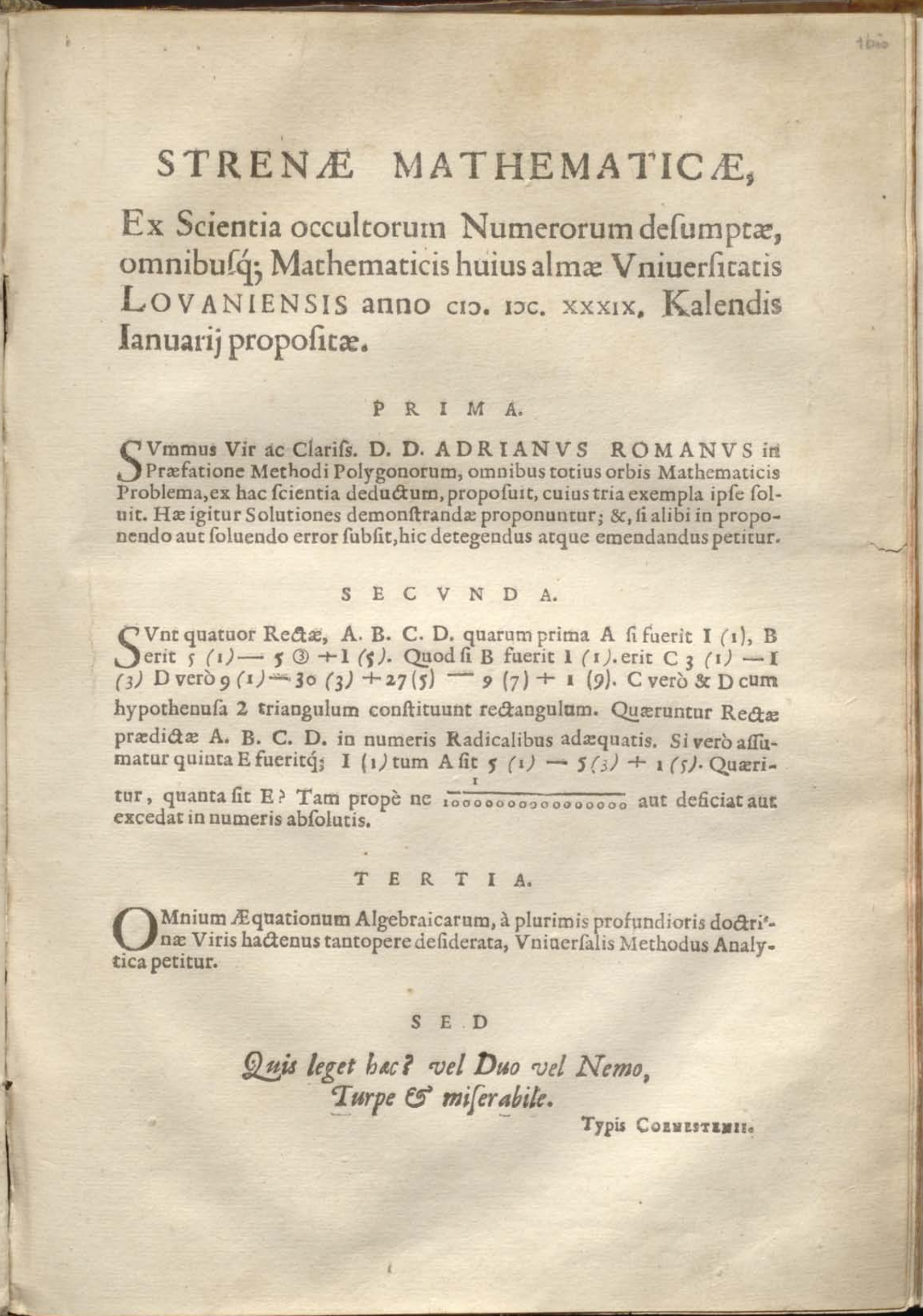}
\caption{Belgian Royal Library KBR, document INC B 635 1 bis}
\label{Romanus-fig2}
\end{figure}    
\newpage

\section{A mathemathical new year's gift}  \label{sec4}
Bockstaele's paper \cite{Bockstaele} describes a mathematical pamphlet that was presented as a new year's gift to all mathematicians of the University of Louvain at the occasion 
of the new year 1639, see Fig. \ref{Romanus-fig2} (courtesy of the Belgian Royal Library). The author\footnote{The author is unknown, but Bockstaele \cite{Bockstaele} names Joannes van der Wegen (a former student of Romanus) or Gerard van Gutschoven (assistant of the Prof. Sturmius, who was teaching mathematics in Louvain at the time) as possible authors.} lists three problems to be solved, but is rather sceptical that anyone will solve them (\textit{maybe two, maybe none}).
The three challenges are
\begin{enumerate}
\item To solve the problem posed by Adrianus Romanus in the preface of his book \textit{Methodi Polygonorum}. Three examples were solved and a proof of these solutions is asked for,
and in case of an error, to find and correct it.
\item To find $A,B,C,D,E$ that satisfy the equations
\begin{eqnarray}
  B &=& 5 A - 5 A^3 + A^5,  \label{eq1} \\
  C &=& 3B - B^3,  \label{eq2} \\
  D &=& 9 B -30 B^3 + 27 B^5 - 9 B^7 + B^9,  \label{eq3} \\
  \lefteqn{C^2+D^2= 4,} \ \ & &  \label{eq4}  \\
  A & = & 5E - 5 E^3 + E^5.  \label{eq5}
\end{eqnarray}   
Express $A,B,C,D$ in radicals and $E$ to a precision of $10^{-16}$.
\item Find a general analytical method to solve an algebraic equation.
\end{enumerate}
For equation \eqref{eq4} the original text says: $C$ and $D$ together with the hypothenusa $2$ form a right-angled triangle. In a picture this means

\begin{figure}[h!]
\centering
\begin{picture}(30,60)(20,-10)
\unitlength=1mm
\put(0,0){\line(1,0){30}}
\put(30,0){\line(0,1){20}}
\put(0,0){\line(3,2){30}}
\put(15,-5){$C$}
\put(35,10){$D$}
\put(13,13){$2$}
\end{picture}
\end{figure}

\noindent so that \eqref{eq4} follows from the Pythagoras theorem.

The problem of Adrianus Romanus in the first challenge is the one we analysed in Section \ref{sec2}--\ref{sec3}. It is clear that the author was aware of the error in Romanus' example 2 and most likely the author knew Vi\`ete's response in \cite{Viete1595}.
The third challenge is over optimistic from the present perspective. Scipione Del Ferro (1465--1526), Niccol\`o Tartaglia (1500--1557), Gerolamo Cardano (1501--1576)  and Lodovico Ferrari (1522--1565) had succeeded in solving cubic and quadratic equations in the sixteenth century, but solving higher order algebraic equations remained an open problem. Nowadays we know that there does not exist a general analytical method to solve algebraic equation of degree five or higher (the Abel-Ruffini theorem), hence the third challenge is an impossible task. However, in 1639 it might have been fully realistic to put this problem forward \cite{Stedall}.

The second challenge is a system of algebraic equations, but again we can recognize Chebyshev polynomials of the first kind:
\begin{eqnarray*}
    2T_3(x/2) &=& x^3 - 3x, \\
    2T_5(x/2) &=& x^5 - 5x^3 + 5x, \\
    2T_9(x/2) &=& x^9 - 9x^7 + 27 x^5 - 30x^3 + 9x.
\end{eqnarray*}
If we use the relation \eqref{Tn-sin} and write the unknowns as
\[   A = 2\sin \theta_1, \quad B=2\sin \theta_2, \quad C = 2\sin\theta_3, \quad D = 2 \sin \theta_4, \quad E = 2 \sin \theta_5, \]
then the equations \eqref{eq1}--\eqref{eq5} become
\begin{eqnarray}
     \sin \theta_2 &=& \sin 5\theta_1, \label{eqt1} \\
     \sin \theta_3 &=& \sin 3\theta_2 , \label{eqt2} \\
     \sin \theta_4 &=& \sin 9\theta_2,  \label{eqt3} \\
     \lefteqn{\sin^2\theta_3 + \sin^2 \theta_4 = 1,} \phantom{\sin \theta_4}  &&  \label{eqt4} \\
     \sin \theta_1 &=& \sin 5\theta_5.   \label{eqt5}
\end{eqnarray}     
From \eqref{eqt4} we find
\[    \cos^2 \theta_3 = \sin^2 \theta_4, \]
so that $\theta_3= \frac{\pi}{2} - \theta_4$ if we look for solutions with $0 < \theta_3,\theta_4 < \pi/2$. 
From \eqref{eqt2} we find $\theta_3 = 3 \theta_2 + 2\pi k$ $(k \in \mathbb{Z})$, and \eqref{eqt3} gives $\theta_4 = 9 \theta_2 + 2\pi \ell$ $(\ell \in \mathbb{Z})$.
Together this gives
\[                3 \theta_2 + 2\pi k  = \frac{\pi}{2} - 9 \theta_2 + 2\pi \ell \quad  \Rightarrow \quad 12 \theta_2 = \frac{\pi}{2} + 2\pi m,  \qquad m \in \mathbb{Z}.  \]
We conclude that
\[     \theta_2 = \frac{\pi}{24} + \frac{2\pi m}{12},  \]
and the smallest positive value for $B$ corresponds to $2\sin(\pi/24)$. The smallest positive value for $C$
then correspond to $C= 2\sin(\pi/8)$ and from \eqref{eq4} we then have $D= 2\cos(\pi/8)$.  From \eqref{eqt1} we also find $\theta_1 = \pi/120$.
In radicals we then have
\begin{eqnarray*}
    A = 2 \sin \frac{\pi}{120} &=& \sqrt{2-\sqrt{2+\sqrt{\frac{3}{16}}+\sqrt{\frac{15}{16}}+\sqrt{\frac{5}{8}-\sqrt{\frac{5}{64}}}}}, \\
    B = 2 \sin \frac{\pi}{24} &=& \sqrt{2-\sqrt{2+\sqrt{3}}},  \\
    C = 2 \sin \frac{\pi}{8} &=& \sqrt{2-\sqrt{2}}, \\
    D = 2 \cos \frac{\pi}{8} &=& \sqrt{2+\sqrt{2}}.  
\end{eqnarray*}      
Finally, for $E$ we have from \eqref{eqt5} that $E = 2\sin \frac{\pi}{600}$.
This cannot be written with radicals since $600=2^3\times 3 \times 5^2$ and hence one needs to solve cubic and quintic equations. In decimals one has
\[   E \approx 0.010471927662839160188 . \]
Clearly we only presented one solution, which corresponds to the geometric meaningful solution where $A,B,C$ and $E$ are sides of a regular polygon,
but there are many more solutions.

\section{Composition identity for Chebyshev polynomials}  \label{sec5}
Chebyshev polynomials of the first kind satisfy a nice composition identity
\begin{equation}  \label{Tmn-TmTn}
    T_n(T_m(x)) = T_{nm}(x), 
\end{equation}    
which follows from \eqref{Tn-cos}. See Rivlin \cite[Chap.~4]{Rivlin}.
This explains why various expressions involving $2\sin (\pi/n)$ or $2\cos (\pi/n)$ can be expressed using a composition of radicals if $n$ is a
nice composite number. The zeros of $T_n$ are given by
\[    \cos \frac{2k-1}{2n} \pi, \qquad   1  \leq k \leq n, \]
so that  $\cos(\pi/2n)$ is one of the $n$ values of $T_n^{-1}(0)$.
If $n=2^j$, then $T_{2^j}(x) = T_2^{\circ j}(x) =T_2(T_2(T_2 \cdots T_2(x)))$, and since $T_2^{-1}(y) = \pm \sqrt{\frac{y+1}{2}}$, the roots
$\cos (2k-1)\pi/2n$ will be a composition of $j$ square roots with a choice of $\pm$ in each step, giving $2^j$ roots. In particular
\[    \cos \frac{\pi}{2^{j+1}} = \sqrt{\frac12 + \frac12 \sqrt{\frac12 + \frac12 \sqrt{\frac12 + \cdots + \frac12 \sqrt{\frac12}}}}, \]
with $j$ square roots. Multiplying by $2$ gives a nicer formula
\[    2\cos \frac{\pi}{2^{j+1}} = \sqrt{2 +  \sqrt{2 + \sqrt{2 + \cdots + \sqrt{2}}}}. \]
For $n=48$ one uses $48=2^4\times 3$ and one of the roots $T_3^{-1}(0) = \sqrt{3}/{2}$ to arrive at
\[   2 \cos \frac{\pi}{48} = \sqrt{2+\sqrt{2+\sqrt{2+\sqrt{3}}}}, \]
as we encountered before.

Recall the solution $\sin(\pi/675)$ of the problem posed by Romanus. Using \eqref{Tn-sin} one has
\[  T_{675}\left(\sin \frac{\pi}{675}\right) = - \sin \pi = 0, \]
so that $\sin (\pi/675)$ is one of the roots $T_{675}^{-1}(0)$. We have that $675 = 3^3 \times 5^2$, hence
\[   T_{675}(x) = T_3(T_3(T_3(T_5(T_5(x))))), \qquad   T_{675}^{-1}(0) = T_5^{-1}(T_5^{-1}(T_3^{-1}(T_3^{-1}(T_3^{-1}(0))))).  \]
The roots $T_3^{-1}(0)$ are $0$ and $\pm \sqrt{3}/2$, but $T_3^{-1}(\sqrt{3}/2)$ involves third roots. This is the reason why $2\sin (\pi/675)$ can not be expressed
using only radicals.
In Section \ref{sec4} the solution for $E$ was found to be $2 \sin (\pi/600)$, which can also be written as $2 \cos (299\pi/600)$. 
From \eqref{Tn-cos} we find that
\[    T_{300}\left( \cos \frac{299\pi}{600} \right) = \cos 299 \frac{\pi}{2} = 0, \]
and hence $\sin(\pi/600)$ is one of the roots $T_{300}^{-1}(0)$. We have
\[    T_{300}(x) = T_2(T_2(T_3(T_5(T_5(x))))), \quad    T_{300}^{-1}(0) = T_5^{-1}(T_5^{-1}(T_3^{-1}(T_2^{-1}(T_2^{-1}(0))))), \]
in which $y=T_2^{-1}(T_2^{-1}(0))$ can be written with radicals, and also $z=T_3^{-1}(y)$ can still be expressed in radicals (even though it involves a cubic equation), but $T_5^{-1}(z)$ requires the solution of a quintic equation and can no longer be expressed in radicals. 

\section{Conclusion}

We gave a few examples where Chebyshev polynomials of the first kind appeared in the 16th (and 17th) century, in work of Fran\c{c}ois Vi\`ete and Adrianus Romanus.
These polynomials appear in the analysis of the sides of regular polygons, which was relevant for approximating $\pi$. The properties that were useful are the
relation with the trigonometric functions \eqref{Tn-cos} and especially \eqref{Tn-sin}, and the composition rule \eqref{Tmn-TmTn} which allows to express
zeros of certain Chebyshev polynomials in radicals using a composition of square roots.
 Other properties of the Chebyshev polynomials
are not to be found, such as the orthogonality relation
\[   \int_{-1}^1 T_n(x)T_m(x) \frac{dx}{\sqrt{1-x^2}} = 0, \qquad m \neq n,  \]
or the differential equation
\[    (1-x^2) T_n''(x) - xT_n'(x) +n^2 T_n(x)=0, \]
which is perfectly normal since differential and integral calculus started much later with Leibniz and Newton.
The three-term recurrence relation is somewhat hidden in Vi\`ete's work \cite{Viete1615}, where he describes how to compute the coefficients of the polynomials.

The problem posed by Adrianus Romanus was to solve the equation $2T_{45}(x/2)=b$ for particular values of $b$. The solution can be obtained by purely geometric
methods, transforming the problem to the comparison of sides of regular polygons. The equation has several (real) solutions, but only the smallest positive solution seemed to
be of interest, since this corresponds to the length of the side of the required polygon. Vi\`ete finds all the positive solutions. The negative solutions are ignored, possibly because
they have no geometric meaning as the length of a side of a polygon\footnote{Other authors gave a geometrical interpretation of negative numbers, such as Bhaskara II (1114--1191) and Albert Girard (1595--1632). There were other more philosophical difficulties towards accepting negative numbers: how can something be less than nothing (0)? I thank the referee for this comment.}.

It is remarkable that one encounters mathematical errors in those texts. Of course there were no referees a the time and printing was a time consuming and manual task.
Nevertheless, the error made by Romanus in his second example is a genuine mathematical error, and surprisingly Vi\`ete's correction was also not the right one.
There are four errors in the table given in Fig.  \ref{figV2} which are likely the result of lack of proofreading\footnote{These errors were corrected in the 1983 translation of T. Richard Witmer \cite[p.~434]{Witmer}, but unfortunately six coefficients which were correct were changed and one coefficient was not changed correctly.} due to the fact that this work was printed in 1615, which is twelve years after Vi\`ete's death in 1603.

Mathematicians liked to challenge each other by posing mathematical problems. The problem posed by Romanus in 1593 was solved fairly quickly by Vi\`ete in 1595, because
Vi\`ete recognized that the problem is equivalent to a problem about regular polygons. Romanus believed that  Ludolf van Ceulen undoubtedly would be able to compute the solution
with many decimals. The new year's gift in Section \ref{sec4} was intended to complain about the quality of mathematics at the University of Louvain where Joannes Storms (Sturmius, 1559--1650) taught mathematics from 1610 until 1646, after Romanus left in 1593. Between 1593 and 1610 there was nobody teaching mathematics, and Sturmius was not a very stimulating
teacher. He was 79 when the pamphlet with the mathematical challenges was distributed and at that time, only two or three students were taking the mathematics classes.

You may get the impression that geometry, algebra and arithmetic were considered as strongly interconnected branches of mathematics in the 16-17th century, but this was certainly not the case. There was a lot of debate and tension about letting arithmetic and algebra enter into the field of geometry \cite[Ch.~6]{Bos}. The merging of algebra and geometry
started with Vi\`ete's \textit{In artem analyticen isagoge} (Introduction to the analytic art)  \cite{Witmer} and culminated in Descartes' \textit{La G\'eometrie} of 1637. 

The main conclusion is that certain algebraic equations of higher degree can still be solved, in particular when they contain Chebyshev polynomials, because one can
transform them to a problem about trigonometric functions and angles and sides of a regular polygon. 
  
 \section*{Acknowledgment}
 The author is grateful to an anonymous referee for some observations, in particular for references \cite{Bos,Cajori,Stedall}

\begin{verbatim}
Walter Van Assche
Department of Mathematics
KU Leuven
Celestijnenlaan 200B box 2400
BE-3001 Leuven
BELGIUM
walter.vanassche@kuleuven.be
\end{verbatim}
\end{document}